\magnification 1250
\pretolerance=500 \tolerance=1000  \brokenpenalty=5000
\mathcode`A="7041 \mathcode`B="7042 \mathcode`C="7043
\mathcode`D="7044 \mathcode`E="7045 \mathcode`F="7046
\mathcode`G="7047 \mathcode`H="7048 \mathcode`I="7049
\mathcode`J="704A \mathcode`K="704B \mathcode`L="704C
\mathcode`M="704D \mathcode`N="704E \mathcode`O="704F
\mathcode`P="7050 \mathcode`Q="7051 \mathcode`R="7052
\mathcode`S="7053 \mathcode`T="7054 \mathcode`U="7055
\mathcode`V="7056 \mathcode`W="7057 \mathcode`X="7058
\mathcode`Y="7059 \mathcode`Z="705A
\def\spacedmath#1{\def\packedmath##1${\bgroup\mathsurround =0pt##1\egroup$}
\mathsurround#1
\everymath={\packedmath}\everydisplay={\mathsurround=0pt}}
\def\nospacedmath{\mathsurround=0pt
\everymath={}\everydisplay={} } \spacedmath{2pt}
\def\qfl#1{\buildrel {#1}\over {\longrightarrow}}
\def\phfl#1#2{\normalbaselines{\baselineskip=0pt
\lineskip=10truept\lineskiplimit=1truept}\nospacedmath\smash {\mathop{\hbox to
8truemm{\rightarrowfill}}
\limits^{\scriptstyle#1}_{\scriptstyle#2}}}
\def\hfl#1#2{\normalbaselines{\baselineskip=0truept
\lineskip=10truept\lineskiplimit=1truept}\nospacedmath\smash{\mathop{\hbox to
12truemm{\rightarrowfill}}\limits^{\scriptstyle#1}_{\scriptstyle#2}}}
\def\diagramme#1{\def\normalbaselines{\baselineskip=0truept
\lineskip=10truept\lineskiplimit=1truept}   \matrix{#1}}
\def\vfl#1#2{\llap{$\scriptstyle#1$}\left\downarrow\vbox to
6truemm{}\right.\rlap{$\scriptstyle#2$}}
\newcount\noteno
\noteno=0
\def\up#1{\raise 1ex\hbox{\sevenrm#1}}
\def\note#1{\global\advance\noteno by1
\footnote{\parindent0.4cm\up{\number\noteno}}{\vtop{\eightpoint\baselineskip12pt\hsize15.5truecm\noindent
#1}}\parindent 0cm}
\def\mono{\lhook\joinrel\mathrel{\longrightarrow}}
\def\longmono{\lhook\joinrel\mathrel{\hfl{}{}}}
\def\iso{\mathrel{\mathop{\kern 0pt\longrightarrow }\limits^{\sim}}}
\font\san=cmssdc10
\def\ext{\hbox{\san \char3}}
\def\sym{\hbox{\san \char83}}

\def\sdir_#1^#2{\mathrel{\mathop{\kern0pt\oplus}\limits_{#1}^{#2}}}
\def\pprod_#1^#2{\raise
2pt \hbox{$\mathrel{\scriptstyle\mathop{\kern0pt\prod}\limits_{#1}^{#2}}$}}

\catcode`\@=11
\font\eightrm=cmr8         \font\eighti=cmmi8
\font\eightsy=cmsy8        \font\eightbf=cmbx8
\font\eighttt=cmtt8        \font\eightit=cmti8
\font\eightsl=cmsl8        \font\sixrm=cmr6
\font\sixi=cmmi6           \font\sixsy=cmsy6
\font\sixbf=cmbx6         \font\tengoth=eufm10       
\font\eightgoth=eufm10 at 8pt  \font\sevengoth=eufm7  
\font\sixgoth=eufm7 at 6pt        \font\fivegoth=eufm5
\newfam\gothfam
\textfont\gothfam=\tengoth \scriptfont\gothfam=\sevengoth
  \scriptscriptfont\gothfam=\fivegoth
\def\goth{\fam\gothfam\tengoth}%
\def\eightpoint{%
  \textfont0=\eightrm \scriptfont0=\sixrm \scriptscriptfont0=\fiverm
  \def\rm{\fam\z@\eightrm}%
  \textfont1=\eighti  \scriptfont1=\sixi  \scriptscriptfont1=\fivei
  \textfont2=\eightsy \scriptfont2=\sixsy \scriptscriptfont2=\fivesy
  \textfont\gothfam=\eightgoth \scriptfont\gothfam=\sixgoth
  \scriptscriptfont\gothfam=\fivegoth
  \def\goth{\fam\gothfam\eightgoth}%
   \textfont\itfam=\eightit
  \def\it{\fam\itfam\eightit}%
  \textfont\slfam=\eightsl
  \def\sl{\fam\slfam\eightsl}%
  \textfont\bffam=\eightbf \scriptfont\bffam=\sixbf
  \scriptscriptfont\bffam=\fivebf
  \def\bf{\fam\bffam\eightbf}%
  \textfont\ttfam=\eighttt
  \def\tt{\fam\ttfam\eighttt}%
  \abovedisplayskip=9pt plus 3pt minus 9pt
  \belowdisplayskip=\abovedisplayskip
  \abovedisplayshortskip=0pt plus 3pt
  \belowdisplayshortskip=3pt plus 3pt 
  \smallskipamount=2pt plus 1pt minus 1pt
  \medskipamount=4pt plus 2pt minus 1pt
  \bigskipamount=9pt plus 3pt minus 3pt
  \normalbaselineskip=9pt
  \setbox\strutbox=\hbox{\vrule height7pt depth2pt width0pt}%
  \let\bigf@nt=\eightrm     \let\smallf@nt=\sixrm
  \normalbaselines\rm}
\catcode`\@=12
\def\pc#1{\tenrm#1\sevenrm}
\def\tx{\kern-1.5pt -}
\def\cqfd{\kern 2truemm\unskip\penalty 500\vrule height 4pt depth 0pt width
4pt\medbreak} 
\def\virg{\raise.4ex\hbox{,}}
\def\decale#1{\smallbreak\hskip 28pt\llap{#1}\kern 5pt}
\def\no{n\up{o}\kern 2pt}
\def\ind{\par\hskip 0.8truecm\relax}
\def\indp{\par\hskip 0.4truecm\relax}
\def\moins{\mathrel{\hbox{\vrule height 3pt depth -2pt width 6pt}}}
\def\rond{\kern 1pt{\scriptstyle\circ}\kern 1pt}
\def\iso{\mathrel{\mathop{\kern 0pt\longrightarrow }\limits^{\sim}}}

\def\Tr{\mathop{\rm Tr}\nolimits}

\input xy
\xyoption{all}
\CompileMatrices
\input amssym.def
\input amssym
\vsize = 25truecm
\hsize = 16truecm
\voffset = -.5truecm
\parindent=0cm
\baselineskip15pt
\overfullrule=0pt
\font\twelvebf=cmbx12
\null\vskip0.3cm
\centerline{\twelvebf Riemannian Holonomy and Algebraic Geometry }
\smallskip
\smallskip \centerline{Arnaud {\pc BEAUVILLE}} \bigskip
\centerline{\eightrm Version 1.1 (25/1/99)}
\vskip1.6cm

{\bf Introduction}
\smallskip

\ind This survey is devoted to a particular instance of the interaction
between Riemannian geometry and algebraic geometry, the study of manifolds with
special holonomy. The holonomy group is one of the most basic objects
associated with a Riemannian metric; roughly, it tells us what are the
geometric objects on the manifold (complex structures, differential
forms, ...) which are {\it parallel} with respect to the metric (see 1.3 for a
precise statement). 
\ind There are two surprising facts about this group. The first one is that,
despite its very general definition, there are  few possibilities -- this is
Berger's theorem (1.2). The second one is that apart from the generic case in
which the holonomy group is $SO(n)$, all other cases appear to be related in
some way  to algebraic geometry.  Indeed the study of {\it compact} manifolds
with special holonomy brings into play some special, and quite interesting,
classes of algebraic varieties: 
 Calabi-Yau, complex symplectic or complex contact manifolds. 
I would like to convince algebraic geometers that this interplay is interesting
on two accounts: on one hand the general theorems on holonomy give deep
results on the geometry of these special varieties; on the other
hand Riemannian geometry provides us with good problems in algebraic geometry
-- see 4.3 for a typical example.

\ind I have tried to make these notes accessible to students with little
knowledge of Riemannian geometry, and a basic knowledge of algebraic geometry. 
Two appendices at the end  recall the basic results of Riemannian
(resp.\ algebraic) geometry which are used in  the text.
\bigskip

{\eightpoint\baselineskip=12pt \leftskip1cm\rightskip1cm\hskip1cm These notes
present a detailed version of the ``Emmy Noether lectures"  I~gave at  Bar
Ilan University (Fall 1998). I want to thank the  Emmy Noether Institute  for
the invitation, and Mina Teicher for her warm hospitality.\par}
\vfill\eject
\centerline{\twelvebf 1. Holonomy}
\medskip
 {\bf 1.1. Definition}
\smallskip

\ind Perhaps the most fundamental object associated to a Riemannian metric on a
manifold $M$ is a canonical connection on the tangent bundle $T(M)$, the {\it
Levi-Civita connection}. A connection gives an isomorphism
between the tangent spaces at infinitesimally near points; more precisely, to each
path $\gamma $ on $M$ with origin $p$ and extremity $q$, the connection associates an
isomorphism $\varphi _\gamma :T_p(M)\rightarrow T_q(M)$ (``parallel transport"),
which is actually an isometry with respect to the scalar products on $T_p(M)$ and
$T_q(M)$ induced by the metric (see App.\ A for more details).  
 If $\delta$ is another path from $q$ to $r$, the
isomorphism associated to the path composed of $\gamma $ and $\delta$ is $\varphi
_\delta\rond\varphi _\gamma $. 
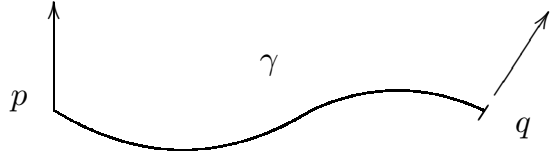
$$\xymatrix @1
{&&&&&&\\
*-{}\ar[u]^<*+<10pt>{p}\ar @/_1pc/@<-2.7pt>
@{-}[rrr]^>>*+<20pt>{\gamma}&&&*-{}\ar
@/^/@<-2.7 pt>
@{-}[rr]&&*-{}\ar@<-2pt>[ur]!<-8pt,-8pt>_<
*+<10pt>{q}|<\hole }$$ 
\ind Let $p\in M$; the above construction associates in particular to every
loop $\gamma $ at $p$ an isometry of $T_p(M)$. The set of all such
isometries is a subgroup $H_p$ of the orthogonal group $O(T_p(M))$, called
the {\it holonomy subgroup} of $M$ at $p$. If $q$ is another point of $M$
and
$\gamma $ a path from $p$ to $q$, we have $H_q=\varphi _\gamma \,
H_p\,\varphi _\gamma ^{-1}$, so that the  $H_p$'s define a unique conjugacy
class $H\i O(n)$; the group $H$ is often  called simply the holonomy group
of $M$. Similarly the representations of the groups $H_p$ on $T_p(M)$ are
isomorphic, so we can talk about the {\it holonomy representation} of $H$.
  \ind There is a variant of this definition, the {\it restricted} holonomy
group, obtained by considering only those loops which are homotopically
trivial. This group actually behaves more nicely: it is a connected, closed
Lie subgroup of $SO(T_p(M))$. To avoid technicalities, we will always assume
that our varieties are {\it simply-connected}, so that the two notions
coincide. We will also usually consider {\it compact} manifolds: this is
somehow the most interesting case, at least for the applications to
algebraic geometry.
\bigskip
{\bf 1.2. The theorems of De Rham and Berger}
\smallskip
\ind With such a  degree of generality  we would expect very few
restrictions, if any, on the holonomy group. This is far from being the
case: thanks to a remarkable theorem of Berger, we can give a complete (and
rather small) list of possible holonomy groups. First of all, let us say
that a Riemannian manifold is {\it irreducible} if its holonomy
representation is  irreducible. 

\smallskip
{\bf Theorem} (De Rham)$.-$ {\it Let $M$ be a compact simply-connected
Riemannian manifold. There exists a canonical decomposition $M\iso
\pprod_{}^{} M_i$, where each $M_i$ is an irreducible Riemannian manifold. Let
$p=(p_i)$ be a point of $M$, and let $H_i\i O(T_{p_i}(M_i))$ be the holonomy
group of $M_i$ at $p_i$; then the holonomy group of $M$ at $p$ is the
product $\pprod_{}^{} H_i$, acting on
$T_p(M)=\pprod_{}^{} T_{p_i}(M_i)$ by the product representation}.
\ind The reader fluent in Riemannian geometry may replace compact by
complete. On the other hand, both completeness and simple connectedness are
essential here. The proof is far from trivial, see for instance [K-N], IV.6.
\bigskip
\ind We are thus reduced to {\it irreducible} (compact, simply-connected)
Riemannian manifolds. Among these are some very classical manifolds, the {\it
symmetric spaces}; they are of the form $G/H$, where $G$ is a compact Lie
group and $H$ is the neutral component  of the fixed locus of an involution
of $G$. These spaces are completely classified, and their geometry  is
well-known; the holonomy group is
$H$ itself. Excluding this case, we get:

\smallskip
{\bf Theorem} (Berger)$.-$ {\it Let $M$ be an irreducible  (simply-connected)
Riemannian manifold, which is not isomorphic to a symmetric space.
Then the holonomy group $H$ of $M$ belongs to the following list:}

\def\tvi{\vrule height 16pt depth 8pt width 0pt}
\def\tv{\tvi\vrule}
\def\q{\quad}
\def\n{\noalign{\hrule}}
\def\h{\hfill}\def\hq{\hfill\quad}
$$\hss\vbox{\offinterlineskip
\halign{\tv$#$&\h\q$#$\hq\tv&  \h\ $#$\h\ \tv  &\h\ #\h\ &\tv#\cr
&H&\dim(M) & {\it metric}&\cr\n
&SO(n)&n & generic &\cr\n
&U(m)&2m & K\"ahler &\cr\n
&SU(m) & 2m & Calabi-Yau &\cr\noalign{\vskip -10pt}
&\scriptstyle{(m\ge 3)}&&&\cr\n
&{\rm Sp}(r)&4r & hyperk\"ahler &\cr\n
&{\rm Sp}(r){\rm Sp}(1)&4r & quaternion-K\"ahler &\cr\noalign{\vskip -10pt}
&\scriptstyle{(r\ge 2)}&&&\cr\n
&G_2&7 &  &\cr\n
&{\rm Spin}(7)&8 &  &\cr\n
}}\hss$$

\smallskip
\ind We have eliminated $SU(2)\ (={\rm Sp}(1))$ and ${\rm Sp}(1){\rm Sp}(1)\
(=SO(4))$ so that a given group appears only once in the list. We should 
point out that a third  exceptional case, ${\rm Spin}(9)\i SO(16)$, appeared
in Berger's list, but  has been eliminated later (see [B-G]).
 \ind Which groups in this list do effectively occur for some compact, 
simply-connected, non-symmetric manifold? That $O(n)$ and $U(m)$ occur is
classical and easy: one starts from an arbitrary Riemannian (resp.\
K\"ahlerian) metric on $M$ and perturbs it in the neighborhood of a point. 
The other groups required much more efforts. The case of  $SU(m)$ is a direct
consequence of   the Calabi conjecture, proved by Yau [Y]; examples with
$H={\rm Sp}(r)$ were found in 1982 [B1], again using Yau's result. Examples
in the last cases, $G_2$ and
${\rm Spin}(7)$, were found only recently [J1, J2]. As for  ${\rm Sp}(1){\rm
Sp}(r)$,  no example is known, and in fact it is generally conjectured that
they should not exist -- we will discuss this in \S 4.

\bigskip

 {\bf 1.3. The holonomy principle}
\smallskip
\ind Before describing  the subgroups which appear in the list,  let us
discuss the geometric meaning of such a restriction on the holonomy. We say
that a tensor field $\theta$ on $M$ is {\it parallel} if for any path
$\gamma $ from $p$ to $q$, the isomorphism $\varphi_\gamma $ transports
$\theta (p)$ onto $\theta (q)$  (this is equivalent to $\nabla^{} \theta=0$, 
see App.\ A). This implies in particular that $\theta (p)$ is invariant
under the holonomy subgroup $H_p$. Conversely, given a tensor $\theta(p)$ on
$T_p(M)$ invariant under $H_p$,  we can transport it at $q$ by any path from
$p$ to
$q$ and obtain a tensor $\theta (q)$  independent of the chosen path; the
tensor field $\theta $ thus constructed is parallel.  We have thus
established: \smallskip {\it Holonomy principle}: Evaluation at $p$
establishes a one-to-one correspondence between parallel tensor fields and
tensors on $T_p(M)$ invariant under $H_p$. \smallskip
\ind In the next sections we will illustrate this principle by going through 
Berger's list. Let us start with the two simplest cases:

\ind a) $H=SO(n)$ means that there are no parallel tensor fields (apart from 
the metric and the orientation). Such a metric is often called {\it generic}.
\smallskip
\ind b) $U(m)$ is the subgroup of $SO(2m)$ preserving a complex structure
$J$ on
${\bf R}^{2m}$ which is orthogonal (that is, $J\in SO(2m),\,J^2=-1$). 
Therefore the manifolds with holonomy contained in $U(m)$ are the Riemannian
manifolds with  a complex structure $J$ which is orthogonal and parallel.
This is one of the classical characterization of {\it K\"ahler} manifolds.
\ind We claimed in the
introduction that compact manifolds with special holonomy are related to
algebraic geometry. In the case of compact K\"ahler manifolds, the link is
provided by the following conjecture:
\smallskip
{\it Is every compact K\"ahler manifold obtained by deformation of a
projective manifold?} 
\smallskip
\ind In dimension 2 this follows from the classification of 
complex surfaces, but nothing is known in dimension $\ge 3$.\smallskip 
\ind We will discuss the groups $SU(m)$, ${\rm Sp}(r)$ and ${\rm Sp}(1){\rm
Sp}(r)$ in the next sections. We will not discuss the exotic holonomies
$G_2$ and ${\rm Spin}(7)$  here; I refer to [J3] for a readable account.

\vskip1cm
\centerline{\twelvebf 2. Calabi-Yau manifolds}
\medskip
\ind We now consider manifolds with holonomy contained in $SU(m)$. We view
$SU(m)$ as the subgroup of $U(m)$ preserving an alternate complex $m$\tx
form on ${\bf C}^m$; therefore a manifold $X$ with holonomy contained
in $SU(m)$ is a K\"ahler manifold (of complex dimension $m$) with a parallel
form of type $(m,0)$. This means that the canonical line bundle $K_X:=\Omega
^m_X$ is flat; in other words, the Ricci curvature (which for a K\"ahler
manifold is just the curvature of  $K_X$) is zero. Thus the
manifolds with holonomy
$SU(m)$ are exactly the  Ricci-flat manifolds. 
\ind It is easy to see that a
parallel form is closed, hence in this case holomorphic: thus the canonical 
bundle $K_X$ of $X$ is  trivial (as a holomorphic bundle). 
Conversely, the Calabi conjecture, proved by Yau [Y], implies that a {\it
Calabi-Yau} manifold, namely a compact, simply-connected
K\"ahler manifold with trivial canonical bundle,  admits a Ricci-flat  metric.
So {\it  the compact (simply-connected) complex manifolds which  admit a 
metric with holonomy contained in $SU(m)$ are the Calabi-Yau manifolds}.
\ind This fact has strong implications in algebraic geometry, in particular
thanks to the following result:
\smallskip
{\bf Proposition} (Bochner's principle)$.-$ {\it On a compact K\"ahler
Ricci-flat manifold, any holomorphic tensor field (covariant or contravariant)
is parallel}. 
\ind The proof rests on the following formula, which follows from a tedious but
straightforward computation ([B-Y], p.\ 142): if $\tau $ is any tensor field,
$$\Delta (\|\tau \|^2)=\|\nabla^{}\tau \|^2\ .$$Therefore $\Delta (\|\tau
\|^2)$ is nonnegative, hence $0$ since its mean value over $X$  is $0$ by
Stokes' formula. It follows that $\tau $ is parallel.\cqfd
\smallskip
\ind As a consequence we get\smallskip 
{\bf Proposition}$.-$ {\it Let $X$ be a compact
K\"ahler manifold, of dimension $m\ge 3$, with holonomy group $SU(m)$.
Then $X$ is projective, and $H^0(X,\Omega ^p_X)=0$ for} $0<p<m$.
\ind {\it Proof}: Let $x\in X$, and $V=T_x(X)$. Using the Bochner and
holonomy principles, we see that the space $H^0(X,\Omega ^p_X)$ can be
identified with the $SU(V)$\tx invariant subspace of $\ext^pV^*$.
Because $SU(V)$ acts irreducibly on $\ext^pV^*$, the invariant subspace is
zero unless $p=0$ or $p=m$. Since $H^0(X,\Omega ^2_X)$ is zero,  $X$ is
projective (App.~B).\cqfd

\smallskip

\ind Manifolds with holonomy ${\rm Sp}(r)$, called {\it hyperk\"ahler}
manifolds, have very special properties; 
we will study them in detail in the next section. Since the only groups in 
Berger's list which are contained in $SU(m)$ are  of the form $SU(p)$ or
${\rm Sp}(q)$,we get the following structure theorem:\smallskip 
{\bf Theorem}$.-$ {\it Any (simply-connected) Calabi-Yau manifold is a
product
$\pprod_i^{}V_i\times \pprod_j^{}X_j$, where:
\indp{\rm a)} Each $V_i$ is a projective Calabi-Yau manifold, with 
$H^0(V_i,\Omega^p_{V_i} )=0$ for $0<p<\dim(V_i)$;
\indp{\rm b)} The manifolds $X_j$  are irreducible hyperk\"ahler}.

\smallskip
\ind  (There is a more general statement for non simply-connected manifolds, see
for instance [B1]).
\bigskip
{\it Further developments}
\ind  Calabi-Yau manifolds have been at the center of a flurry of activity
in the last 10 years, principally under the influence of mathematical
Physics. The key word here is {\it mirror symmetry}, a (conjectural)
duality  between families of  Calabi-Yau manifolds. I will not try to be
more precise, because this goes far beyond the scope of these notes. An
excellent reference is the booklet [V]. The current trend puts the emphasis
on the symplectic, rather than algebro-geometric, aspect [S-Y-Z].

\vskip1cm
\centerline{\twelvebf 3. Symplectic manifolds}
\medskip

 {\bf 3.1. Hyperk\"ahler versus symplectic}\smallskip 
\ind The group ${\rm
Sp}(r)$ is the quaternionic unitary group, that is, the
group of ${\bf H}$\tx linear automorphisms of ${\bf H}^r$ which preserve the
standard hermitian  form $\psi (z,z')=\sum_r z_i\bar z'_i$. Viewing  ${\bf
H}^r$ as ${\bf R}^{4r}$ realizes ${\rm Sp}(r)$ as a subgroup of the
orthogonal group ${\rm SO}(4r)$. The manifolds of dimension $4r$ with
holonomy ${\rm Sp}(r)$ are called {\it hyperk\"ahler} manifolds.
\ind There are two ways of making this definition explicit.  We can
characterize ${\rm Sp}(r)$ as the subgroup of orthogonal transformations of 
${\bf R}^{4r}$  which are linear with respect to the complex structures
$I,J,K$ (here $(1,I,J,K)$ is the standard basis of ${\bf H}$ over ${\bf R}$,
with $IJ=-JI=K$). By the holonomy principle, hyperk\"ahler manifolds are
therefore characterized by the existence of $3$ complex structures $I,J,K$,
with $IJ=-JI=K$, such that the metric is K\"ahler with respect to each of
these. Actually any pure quaternion $aI+bJ+cK$ with $a^2+b^2+c^2=1$ defines
such a structure, so hyperk\"ahler manifolds admit a family of complex
K\"ahler structures parametrized by the sphere ${\bf S}^2$ (hence their
name).  
\ind A second way to look at ${\rm Sp}(r)$ is to give a special role to  one
of these complex structures, say $I$, and to view ${\bf H}$ as ${\bf C}(J)$
(and ${\bf C}$ as ${\bf R}(I)$). We identify ${\bf H}^r$ with ${\bf
C}^r\oplus {\bf C}^rJ={\bf C}^{2r}$. The hermitian form $\psi $ can be
written as $h+\varphi J$, where $h$ is the standard (complex) hermitian
form  and $\varphi $ the standard  ${\bf C}$\tx bilinear symplectic form on
${\bf C}^{2r}$. Therefore ${\rm Sp}(r)$ is the intersection in $SO(4r)$ of
the unitary group $U(2r)$ and the complex symplectic group ${\rm Sp}(2r,{\bf
C})$ (incidentally, this implies that ${\rm Sp}(r)$ is a maximal compact
subgroup of ${\rm Sp}(2r,{\bf C})$, which is the reason for the notation). 
\ind In terms of holonomy, this means that once a preferred complex
structure has been chosen, a hyperk\"ahler manifold can be characterized as a
K\"ahler manifold with a parallel non-degenerate 2-form of type $(2,0)$. As
above this 2-form must be holomorphic, hence it is  a (complex) {\it
symplectic structure}, that is a closed\note{The closedness condition is
automatic for compact K\"ahler manifolds.}, holomorphic, everywhere
non-degenerate 2-form. Conversely, let $X$ be a compact K\"ahler manifold of
(complex) dimension $2r$, with a complex  symplectic structure $\varphi$;
then $X$ is a Calabi-Yau manifold (because $\varphi^r$ does not vanish),
hence admits a Ricci-flat metric, for which the form $\varphi $ is parallel.
If moreover we require the holomorphic 2-form $\varphi$ to be unique up to a
scalar, the holonomy of $X$ is exactly ${\rm Sp}(r)$.
We will call such a manifold {\it K\"ahler symplectic}, to emphasize that
we have chosen a particular complex structure.
\bigskip
{\bf 3.2. The two standard series}
\smallskip
\ind A typical example of a K\"ahler symplectic manifold  is a K3
surface, that is a compact (simply-connected) K\"ahler surface with trivial
canonical bundle. Note that in the statement of Berger's theorem I have
deliberately chosen to view the group $SU(2)$ as symplectic $(={\rm
Sp}(1))$ rather than unitary: we  will see that the theory of K3 surfaces
is an accurate model for the study of complex symplectic manifolds.   For
a long time no other example  has been known, and it was even conjectured 
that such manifolds should not exist (see [Bo1]). In 1982 Fujiki gave an
example in dimension 4, which I generalized in any dimension -- in fact I
constructed two series of examples [B1]. Let me explain these examples.
\def\pr{\mathop{\rm pr}\nolimits}
\ind Start from a K3 surface $S$, with a holomorphic nonzero 2-form $\varphi 
$.  The product $S^r$ admits a natural symplectic form, namely
$\pr_1^*\varphi +\ldots +\pr_r^*\varphi $; but there are others, since we
may  take as well  any expression $\lambda _1\pr_1^*\varphi+\ldots+\lambda
_r\pr_r^*\varphi$ with
$\lambda _1,\ldots,\lambda _r$ in ${\bf C}^*$. A natural way to eliminate 
those is to ask for ${\goth S}_r$\tx invariant 2-forms, which amounts to
consider instead of $S^r$ the symmetric product $S^{(r)}:=S^r/{\goth S}_r$.
\ind Unfortunately this quotient is singular as soon as $r$ is greater than 
1; but it admits a nice desingularization, the {\it Douady space} 
$S^{[r]}$ which para\-meterizes the finite subspaces of $S$
of length $r$ (when $S$ is projective this is known as 
 the Hilbert scheme). We can view $S^{(r)}$ as the space of finite subsets
$E\i S$ with a positive multiplicity $m(p)$ assigned to each point $p$ of 
$E$, in such  a way that  $\sum_{p\in E}m(p)=r$. The natural map
$\varepsilon : S^{[r]}\rightarrow S^{(r)}$ which associates to a subspace
$Z$ of $S$ its set of points counted with multiplicity turns out to be
holomorphic; it induces an isomorphism  on the open subset $S^{[r]}_{\rm o}$
of $S^{[r]}$ parameterizing those subspaces which consist of $r$ distinct
points. 
\ind It is then easy to show that the 2-form $\pr_1^*\varphi
+\ldots +\pr_r^*\varphi $, which lives naturally on $S^{[r]}_{\rm o}$, 
extends to a symplectic form on  $S^{[r]}$, unique up to a scalar, and that
$S^{[r]}$ is simply-connected. Moreover $S^{[r]}$ is  K\"ahler as a
consequence of a general result of Varouchas [Va]. In other words, {\it the
Douady space $S^{[r]}$ is a
$(2r)$\tx dimensional irreducible symplectic manifold}.
\ind We can perform the same construction starting from a 2-dimensional 
complex torus
$T$: the Douady space $T^{[r]}$ is again symplectic, however it is not
simply-connected. In fact it admits a smooth surjective map
$S:T^{[r]}\rightarrow T$, which is the composite of $\varepsilon
:T^{[r]}\rightarrow T^{(r)}$ and of the sum map
$T^{(r)}\rightarrow T$. The fibre $K_{r-1}=S^{-1}(0)$ is a
simply-connected, irreducible symplectic manifold of dimension $2r-2$.
\smallskip 
\ind Thus we get two series of examples in each dimension. The first thing
to  look at, for algebraic geometers, is  their deformations: there are
some obvious ones obtained by deforming the surface $S$ (or $T$), but it 
turns out that we get more than those. In fact, in the moduli space
parameterizing all deformations of the manifolds we found, those of the
form  $S^{[r]}$ for some
$K3$ surface $S$ form a {\it hypersurface}, and similarly for $K_r$ (this is, 
of course, for $r\ge 2$).	
\ind This is seen as follows. First of all, the universal deformation space
of a symplectic manifold $X$ is smooth, of dimension $\dim H^1(X,T_X)$. This
is a general result for Calabi-Yau manifold, due to Tian and Todorov (see
[T]); in the particular case of  symplectic manifolds it had been proved
earlier by Bogomolov [Bo1]. Since $X$ is symplectic,
 the tangent sheaf
$T_X$ is isomorphic to $\Omega^1_X$, hence
$$\dim H^1(X,T_X)=\dim H^1(X,\Omega^1 _X)=b_2(X)-2\ .$$
An easy computation gives $b_2(S^{[r]})=b_2(S)+1$ and $b_2(K_r)=b_2(T)+1$ for
$r\ge 2$, hence our assertion .
\ind We will say that a symplectic
manifold is {\it of type} $S^{[r]}$, or $K_r$, if it can be obtained by
deformation of $S^{[r]}$, or $K_r$. As an example, we proved in [B-D]
that the variety of lines contained in a smooth cubic hypersurface $V$ of
${\bf P}^5$ is of type $S^{[2]}$, but it is not isomorphic to $S^{[2]}$ if
$V$ is general enough.\bigskip
{\bf 3.3. Other examples}\smallskip 
\ind  Shortly after the two series were discovered, Mukai showed that they
fit into an elegant  construction which looks much more general [M]. He
proved that the moduli space of stable vector bundles on a K3 or abelian
surface  $S$, with fixed rank and Chern classes, is smooth and admits a
symplectic form. The idea is quite simple. The smoothness follows from a
standard  obstruction argument: one shows that the obstructions to deform
$E$ infinitesimally are the same as the obstructions to deform $\det E$,
which vanish. Now the tangent space to the moduli space at $E$ is $H^1(S,
{\cal E}nd(E))$, and the symmetric form $(u,v)\mapsto \Tr uv$ on ${\cal
E}nd(E)$ gives rise to a skew-symmetric pairing
$$H^1(S, {\cal E}nd(E))\otimes H^1(S, {\cal E}nd(E))\longrightarrow
H^2(S,{\cal O}_S)\cong {\bf C}$$
which is non-degenerate by Serre duality, and provides the required symplectic
form.\smallskip 
\ind If we want to exploit this construction to give new examples of 
symplectic manifolds, we need to fulfill the following requirements:
\indp a) Our moduli space $M$ should be compact. This is achieved by including
in $M$   {\it stable sheaves}, and choosing the polarization so that all
semi-stable sheaves are actually stable. I refer for instance to [H-L] for the
details. 
\indp b) $M$ should be simply-connected, and satisfy $\dim H^0(M,\Omega
^2_M)=1$.  This was proved  in [OG1]. Observe that both properties
are invariant by deformation, and also under birational equivalence. 
O'Grady deforms $S$ to a  surface $S_e$ admitting an elliptic pencil, with a
suitable polarization; then a detailed analysis shows that the moduli space is
birational to $S_e^{[r]}$ for some $r$.
\ind 
So $M$  is a  symplectic manifold, and more precisely  a
deformation of a symplectic manifold  of type  $S^{[r]}$. This is actually 
more than we would wish: Huybrechts proved recently that two birational
symplectic manifolds are deformation of each other -- we will discuss this in
detail in 3.5. Therefore {\it the moduli space $M$ is of type} $S^{[r]}$, and
thus does not provide any new example.
\ind  When Huybrechts' result appeared, it implied that all known examples
of K\"ahler symplectic manifolds were of type $S^{[r]}$ or
$K_r$. Since then a new example has been constructed by O'Grady [OG2], of
dimension 10, by desingularizing a  {\it singular} moduli space of vector
bundles on a K3. It still remains an intriguing and very interesting problem
to construct more examples. As we will see in the next sections, we know a
lot about the geometry of K\"ahler symplectic manifolds; it is somewhat
embarrassing to have so few examples.

 \bigskip {\bf 3.4. The
period map} \smallskip
\ind For K3 surfaces the theory of the period map gives us a fairly complete
picture of the moduli space, thanks to the work of Shafarevich and
Piatetski-Shapiro, Burns and Rapoport, Todorov, Looijenga, Siu -- I refer to
[B2] for a survey. The idea is to encode a K3 surface $S$ by its Hodge
decomposition (see App.\ B)$$H^2(S,{\bf C})=H^{2,0}\oplus H^{1,1}\oplus
H^{0,2}\ ,$$ which is
determined by the position of the line $H^{2,0}$ in $H^2(S,{\bf C})$
(we have $H^{0,2}=\overline{H^{2,0}}$, and $H^{1,1}$ is the orthogonal of 
$H^{2,0}\oplus  H^{0,2}$ for the intersection product). The point is that 
$H^2(S,{\bf C})$ depends only on the topology of $S$, while $H^{2,0}$ depends
heavily on the complex structure: we have $H^{2,0}={\bf C}\varphi $, where
$\varphi $ is the De Rham class of a non-zero
holomorphic 2-form on $S$ (unique up to a constant).
\ind To be more precise, we denote
by $L$ a lattice isomorphic to $H^2(S,{\bf Z})$ (this is the unique even
unimodular lattice of signature $(3,19)$, but we will not need this). A {\it
marked} K3 surface is a pair $(S,\sigma )$ of a K3 $S$ and a lattice
isomorphism $\sigma :H^2(S,{\bf Z})\rightarrow L$. The first (easy) result is
that there is an analytic manifold ${\cal M}_L$ which is a {\it fine moduli
space} for marked K3's: that is, there is a universal family $u:{\cal
U}\rightarrow {\cal M}_L$ of marked K3's over ${\cal M}_L$, such that any family 
${\cal S}\rightarrow T$ of marked K3's is the pull-back of
 $u$ through a classifying map $T\rightarrow {\cal M}_L$. Note however that
${\cal M}_L$ is {\it not Hausdorff} -- a rather surprising fact that we will
explain later (3.5).
\ind The advantage of working with ${\cal M}_L$ is that we can now compare the
Hodge structures of different surfaces. Given $(S,\sigma )$, we extend $\sigma
$ to an isomorphism $H^2(S,{\bf C})\rightarrow L_{\bf C}$ and put\note{
We denote as usual by ${\bf P}(V)$ the space of lines in a vector space $V$, and
 by $[v]\in {\bf P}(V)$ the line spanned by a nonzero vector   $v$ of $ V$.}
$$\wp(S,\sigma )=\sigma(H^{2,0})=\sigma ([\varphi]) \in {\bf P}(L_{\bf
C})\ .$$
 The map
$\wp$ is called the {\it period map}, for the following reason: choose a basis
$(e_1,\ldots,e_{22})$ of $L^*$, so that $L_{\bf C}={\bf C}^{22}$. Put $\gamma
_i={}^t\sigma(e_i)$, viewed as an element of 
$H_2(S,{\bf Z})$; then $$\wp(S,\sigma )=\Bigl(
\int_{\gamma _1}\varphi\ :\ \ldots\ :\  \int_{\gamma _{22}}\varphi
\Bigr)\ \in {\bf P}^{21}\ ;$$ the numbers $\int_{\gamma _i}\varphi $ are
classically called the ``periods" of $\varphi $.  \ind
Since $\varphi $ is holomorphic we have $\varphi \wedge\varphi =0$ and
$\int_S\varphi \wedge\bar\varphi >0$. In other words, $\wp(S,\sigma )$ lies in
the subvariety $\Omega_L$ of ${\bf P}(L_{\bf C})$, called the {\it period
domain}, defined by  $$\Omega_L =\{[x]\in {\bf P}(L_{\bf C})\ |\ x^2=0\ ,\
x\bar x>0  \}\ .$$

{\bf Theorem}$.-$ 1) $\wp:{\cal M}_L\rightarrow \Omega_L$ {\it is \'etale
and surjective}.
\ind 2) {\it If $\wp(S,\sigma )=\wp(S',\sigma' )$, the surfaces $S$ and $S'$
are isomorphic}.
\ind Note that this does {\it not} say that $\wp$ is an isomorphism (otherwise
${\cal M}_L$ would be Hausdorff!): the same K3 with different markings can have
the same period. There is a more precise statement which describes exactly the
fibres of $\wp$ (see for instance [P], p.~142, prop.~2).\smallskip 
{\bf Corollary}$.-$ {\it Every $K3$ surface is a deformation of a projective
one}.
\ind {\it Proof}: Write $\varphi=\alpha+i\beta$, with $\alpha,\beta\in
H^2(S,{\bf R})$. The condition $[\varphi]\in \Omega_L$ translates as
$\alpha^2=\beta^2>0$, $\alpha.\beta=0$. It follows that the classes $[\varphi]$
with $\alpha,\beta\in H^2(S,{\bf Q})$ are dense in $\Omega_L$. The corresponding
surfaces are dense in ${\cal M}_L$; they have $H^{1,1}=({\bf C}\alpha\oplus{\bf
C}\beta)^\perp$ defined over ${\bf Q}$, hence they are projective (App.\ B).\cqfd
\ind Note that we only need an easy part of the theorem, namely the fact that
$\wp$ is \'etale. 
\medskip
 
\ind We want to apply the same approach for any K\"ahler symplectic manifold $X$. We
still have the Hodge decomposition
$$H^2(X,{\bf C})=H^{2,0}\oplus H^{1,1}\oplus H^{0,2}\quad\hbox{with }\ H^{2,0}=
{\bf C}\varphi\ .$$
What seems to be lacking is the quadratic form, but in fact it is still there: 
I showed in [B1] that the point $[\varphi]\in{\bf P}(H^2(X,{\bf C}))$ must lie
in a hyperquadric, which is rational over ${\bf Q}$; this implies that there
exists a canonical quadratic form $q:H^2(X,{\bf Z})\rightarrow {\bf Z}$. It has
the following properties (see [B1] and [H1]):
\indp {\it a})  $q$ is non-divisible, non-degenerate, of signature $(3,b_2-3)$;
\indp {\it b}) there exists a positive integer $d_X$ such that $\alpha
^{2r}=d_X\,q(\alpha )^r$ for all $\alpha \in H^2(X,{\bf Z})$;
 \indp {\it c}) $q(\varphi )=0$, and $q(\varphi +\bar \varphi )>0$. 
 \ind We can now mimic the K3 case. Let
$L$ be a lattice; we define as before the moduli space ${\cal M}_L$ of pairs
$(X,\sigma )$, where $X$ is K\"ahler symplectic manifold and 
$\sigma :H^2(X,{\bf Z})\rightarrow L$ a lattice isomorphism. We still have a
natural structure of analytic (non-Hausdorff) manifold on ${\cal M}_L$ (it is
however no longer a fine moduli space in general).  To each element
$(X,\sigma)$ of ${\cal M}_L$ we associate $$\wp(X,\sigma
)=\sigma(H^{2,0})=\sigma ([\varphi ])
\in {\bf P}(L_{\bf C})\ .$$
As above, if we choose a basis $(e_1,\ldots,e_{b})$ of $L^*$, the element 
$\wp(X,\sigma )$ is given by the ``periods" 
$\int_{\gamma _i}\varphi$, with $ \gamma_i={}^t\sigma (e_i)$.
   \ind By property
{\it c}) of $q$,  $\wp(X,\sigma )$ lies in
the subvariety $\Omega_L$ of ${\bf P}(L_{\bf C})$  defined by 
$$\Omega_L =\{[x]\in {\bf P}(L_{\bf C})\ |\ q(x)=0\ ,\ q(x+\bar x)>0  \}\ .$$
{\bf Theorem}$.-$  $\wp:{\cal M}_L\rightarrow \Omega_L$ {\it is \'etale
and surjective}.
\ind The fact that $\wp$ is \'etale follows from the (easy) computation of
its tangent map. The  much more delicate surjectivity has been proved by
Huybrechts [H1].\smallskip 
\ind Using the easy part of the theorem and  the same
argument as for K3 surfaces we obtain:\par
{\bf Corollary}$.-$ {\it Every K\"ahler symplectic manifold is a deformation of 
a projective one}.\medskip
\ind On the other hand, the Torelli problem is still wide open.
There are examples, due to Debarre [De], of nonisomorphic K\"ahler symplectic
manifolds with the same periods; the best one can hope for is:\smallskip 
{\bf Torelli problem}$.-$ {\it If $\wp(X,\sigma )=\wp(X',\sigma' )$, are $X$ and $X'$
birational}?

\bigskip
{\bf 3.4. Birational symplectic manifolds}
\smallskip
\ind The fact that the moduli space ${\cal M}_L$ of marked K3 surfaces is
non-Hausdorff goes back to a famous example of Atiyah [A]. Start with a family
$f:{\cal X}\rightarrow D$ of K3 surfaces over the unit disk, such that the
total space ${\cal X}$ is smooth, the surface ${\cal X}_t$ is smooth for
$t\not=0$ and ${\cal X}_0$ has an ordinary double point $s$: near $s$ we can
find local coordinates $(x,y,z)$ such that $f(x,y,z)=x^2+y^2+z^2$. Pull back
$f$ by the covering $t\mapsto t^2$ of the disk: we obtain a new family ${\cal
Y}\rightarrow D$, where now ${\cal Y}$ has an ordinary double point
$x^2+y^2+z^2=t^2$. Blowing up $s$ in ${\cal Y}$ we get a smooth threefold
$\widehat{\cal Y}$ with a smooth quadric $Q$ as exceptional divisor; we can
now blow down $Q$ along each of its two rulings to get smooth  threefolds
${\cal Y}'$, ${\cal Y}''$, which are {\it small resolutions} of ${\cal Y}$:
the singular point $s$ has been blown-up to a line. 
 
$$
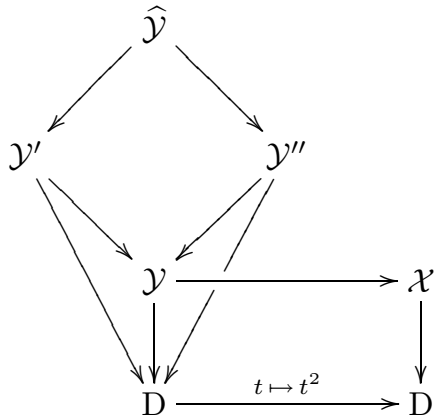
\xymatrix{
& {\widehat{\cal Y}} \ar[dl]\ar[dr] & &  \\
{\cal Y}'\ar[dr]\ar[ddr] & & {\cal Y}''\ar[dl]
\ar[ddl]  |! { [dl]; [dr]}\hole
&
\\ &{\cal Y}\ar[rr]
\ar[d] & & {\cal X}\ar[d]\\
& D\ar[rr]^{t\,\mapsto\, t^2} & & D\\}$$ \smallskip 
 \ind The two fibrations ${\cal Y}'\rightarrow D$ and ${\cal Y}''\rightarrow D$
are smooth; their fibres at $0$ are both isomorphic to the blow up of ${\cal X}_0$ at
$s$. By construction they coincide  above $D\moins\{0\}$, but it is easily
checked that the isomorphism does not extend over $D$.
\ind The local systems  $H^2({\cal Y}'_t, {\bf Z})_{t\in D}$ and  $H^2({\cal Y}''_t,
{\bf Z})_{t\in D}$ are constant, and coincide over $D\moins\{0\}$; choosing
compatible trivializations  we get two non-isomorphic families of marked K3 surfaces
on $D$, which coincide on $D\moins\{0\}$. The corresponding maps $D\rightarrow {\cal
M}_L$ coincide on $D\moins\{0\}$, but take different values at $0$. In other
words, the marked surfaces ${\cal Y}'_0$ and ${\cal Y}''_0$ give {\it non-separated}
points in the moduli space ${\cal M}_L$ (every neighborhood of one of these
points contains the other one). 
\medskip
\ind To explain the analogous construction for higher-dimensional
symplectic manifolds, let us first describe, in the simplest possible case, the
{\it elementary transformations} discovered by Mukai [M]. We 
start with a symplectic manifold $X$, of dimension
$2r$, containing a submanifold $P$ isomorphic to ${\bf P}^r$. The 2-form
$\varphi$ restricted to $P$ vanishes (in fancy words, $P$ is a Lagrangian
submanifold); therefore we have a commutative diagram of exact
sequences\vskip-10pt
$$\diagramme{
0\ \longrightarrow & T_P & \hfl{}{} & T_{X\,|P}&  \hfl{}{} & N_{P/X}
& \longrightarrow\  0\cr & \vfl{}{} & & \vfl{\varphi}{} & & \vfl{}{} &\cr
0\ \longrightarrow & N_{P/X}^* & \hfl{}{} & \Omega^1_{X\,|P}&  \hfl{}{} & \Omega^1_{P} &
\longrightarrow\  0\cr }$$\vskip-5pt
in which all vertical arrows are isomorphisms. In particular, $N_{P/X}$ is
isomorphic to $\Omega^1_P$.
\ind Now blow-up $P$ in $X$:\vskip-20pt
$$\diagramme{
E &\longmono &\widehat{X}\cr
 \vfl{}{} & &  \vfl{}{} \cr
P &\longmono &X \cr
}$$ \vskip-10pt
The exceptional divisor $E$ is by definition the projective normal bundle\note{We
use the standard differential-geometric notation:  if
$F$ is a vector bundle on a variety $B$, we put ${\bf P}(F)=\cup_{b\in B}{\bf
P}(F_b)$ (see footnote $^2$).}
$\ {\bf P}(N_{P/X})$, which by the above remark is
isomorphic to the projective cotangent bundle ${\bf P}T^*(P)$; thus we can view $E$
as the  variety 
 of pairs $(p,h)$ with $p\in P$, $h\in P^*$ (the  space of hyperplanes in
$P$) and $p\in h$. This is clearly symmetric:  $E$ is also isomorphic to
${\bf P}T^*(P^*)$, and in fact, using a classical contractibility criterion (due to
Fujiki and Nakano in this context), we can blow down $E$ onto
$P^*$ and get a new symplectic manifold $X'$, called the elementary transform of
$X$ along
$P$.  The map $X\dasharrow X'$ is a typical example of a birational map between
symplectic manifolds which is not an isomorphism. Note that it is not known
whether $X'$ is always K\"ahler.
\smallskip 
\ind Now suppose we deform $X$ in a family ${\cal X}\rightarrow D$.  We have
an exact sequence of normal bundles
$$0\rightarrow N_{P/X}\cong \Omega^1_P\longrightarrow N_{P/{\cal
X}}\longrightarrow N_{X/{\cal X}}\cong{\cal O}_P\rightarrow 0\ .$$
The class of this extension lives in $H^1(P,\Omega^1_P)$; a straightforward
computation shows that it is the restriction of the tangent vector in the
deformation space of $X$ provided by the deformation ${\cal X}\rightarrow D$
(remember that this  tangent vector belongs to
$H^1(X,T_X)\cong H^1(X,\Omega^1_X)$). Choose ${\cal X}$ so that  this
tangent vector does not vanish on $P$, for instance is a K\"ahler class in
$H^1(X,\Omega^1_X)$. Then the above extension is the non-trivial Euler
extension
$$0\rightarrow \Omega^1_P\longrightarrow V^*\otimes_{\bf C} {\cal
O}_P(-1)\longrightarrow N_{X/{\cal X}}={\cal O}_P\rightarrow 0\ ,$$
where $P={\bf P}(V)$. 
So we get an isomorphism $N_{P/{\cal
X}}\cong V^*\otimes_{\bf C}{\cal O}_P(-1)$. Thus if we blow-up $P$ in
${\cal X}$, the exceptional divisor ${\cal E}$ is isomorphic to $P\times P^*$.
As before we can blow-down ${\cal E}$ onto $P^*$ and get a  manifold
${\cal X}'$ with a smooth map ${\cal X}'\rightarrow D$, whose fibre at $0$ is
isomorphic to $X'$. Again the two families coincide above $D\moins\{0\}$. 
Therefore if  $X'$ is K\"ahler, $X$ {\it and} $X'$
(with appropriate markings) {\it give non-separated points in the moduli 
space} ${\cal M}_L$.
\smallskip 	
\ind This example, due to D. Huybrechts, was the point of departure of his
investigation of birational symplectic manifolds. The outcome
is:\smallskip   {\bf Theorem} (Huybrechts, [H1, H2])$.-$ {\it Let
$X$, $X'$ be two birational K\"ahler symplectic manifolds. There
exists smooth families ${\cal X}\rightarrow D$ and
${\cal X}'\rightarrow D$ which are isomorphic over
$D\moins\{0\}$ and such that ${\cal X}_0 $ is  isomorphic to
$X$ and ${\cal X}'_0 $ to
 $X'$.}
	\ind As before it follows that $X$  and $X'$,
with appropriate markings, give non-separated points in the moduli
space ${\cal M}_L$.
 Conversely, Huybrechts also proves that if $(X,\sigma)$ and
$(X',\sigma')$ are non-separated points in ${\cal M}_L$, the 
 manifolds $X$ and $X'$ are birational [H1].\smallskip 
{\bf Corollary}$.-$ {\it Two  K\"ahler
symplectic manifolds which are birational are diffeomorphic}.
\ind It is interesting to compare this statement with the following result 
of Batyrev [Ba]:
\smallskip 
{\bf Proposition}$.-$  {\it Two Calabi-Yau manifolds which are birational
have the same Betti numbers}.
\ind The proof is (of course) completely different: it proceeds by reduction 
to characteristic $p$. Note that the two Calabi-Yau manifolds need not be
diffeomorphic, as shown by an example of Tian and Yau (see [F], example 7.7).

\bigskip

{\bf 3.5. Further developments}
\ind  K\"ahler
symplectic manifolds have been much studied in the recent years; there are two
directions which I would like to emphasize. The structure of the cohomology
algebra has been  studied by Verbitsky; we will follow the elegant presentation
of  Bogomolov [Bo2]. 
\smallskip  {\bf
Proposition}$.-$ {\it Let
$X$ be a K\"ahler symplectic manifold of dimension $2r$, and let $A$ be the
subalgebra  of
$H^*(X,{\bf Q})$ spanned by $H^2(X,{\bf Q})$. Then $H^*(X,{\bf Q})=A\oplus
A^{\perp}$, and $A$ is the quotient of $\sym^* H^2(X,{\bf Q})$ by the ideal
spanned  by the elements
$x^{r+1}$ for all  
$x\in H^2(X,{\bf Q})$ with
$q(x)=0$. }
\ind Let $Q$ be the quadric $q(x)=0$ in $H^2(X,{\bf C})$. Since the period map is
\'etale (3.3), we know that there is an  open subset $V$ of $Q$ such that
any element
 of $V$ is the class of a 2-form on $X$, holomorphic with respect to some
complex structure on  $X$. This implies $x^{r+1}=0$ for $x\in V$, and therefore
for all $x\in Q$ by analytic continuation.
\ind The rest of the proof is purely algebraic. Given a vector space $H$ over 
${\bf Q}$ with a non-degenerate quadratic form $q$, we consider the algebra 
$A_r(H,q)$ quotient of $\sym^* H$ by the ideal spanned by the elements
$x^{r+1}$ for all $x\in H$ with $q(x)=0$.  Using the representation theory of
$O(H,q)$, one proves that $A_r(H,q)$ is a {\it Gorenstein} algebra; more precisely
$A^{2r}_r(H,q)$ is one-dimensional, and the pairing $A^i_r(H,q)\times
A^{2r-i}_r(H,q)\rightarrow A^{2r}_r(H,q)\cong{\bf Q}$ is non-degenerate for
each $i$. 
\ind Put $H=H^2(X,{\bf Q})$. By the geometric property  above we get a
 ring homomorphism
$ A_r(H,q)\rightarrow  H^*(X,{\bf Q})$. Its kernel is an
ideal of $A_r(H,q)$; if it is non-zero, it  contains the minimum ideal
$A^{2r}_r(H,q)$, so the map $\sym^{2r}H\rightarrow H^{4r}(X,{\bf Q})$ is
zero -- which is impossible since  $\omega^{2r}\not=0$ for a K\"ahler class
$\omega$. Hence $A$ is isomorphic to $ A_r(H,q)$; since the restriction of
the intersection form on  $H^*(X,{\bf Q})$ to $A$ is non-degenerate, we have 
$H^*(X,{\bf Q})=A\oplus A^{\perp}$.\cqfd
\medskip
\ind Another exciting recent development is the construction by Rozansky
and Witten of invariants of 3-manifolds associated to any compact
hyperk\"ahler manifold ([R-W]; an account more readable for an algebraic
geometer appears in [K]). By the advanced  technology of 3-dimensional
topology, defining such invariants amounts to associate a complex number 
(a ``weight")  to
each trivalent graph, in
such a way that a certain identity, the so-called IHX relation, is satisfied.
The weights associated by Rozansky
and Witten to a hyperk\"ahler manifold are sort of  generalized Chern
numbers, which certainly deserve further study. Some explicit computations
have been done by Hitchin and Sawon (to appear). 
\vskip1cm
\centerline{\twelvebf 4. Quaternion-K\"ahler manifolds}\medskip
{\bf 4.1. The twistor space}
\smallskip 
\ind The group ${\rm Sp}(1)$ is the group of quaternions of norm $1$; it acts on  ${\bf
H}^r$  by homotheties.  Since
${\bf H}$ is not commutative, it is {\it not} contained  in the unitary
group ${\rm Sp}(r)$, but it of course commutes with ${\rm Sp}(r)$. A manifold of 
dimension $4r$ is said to be {\it quaternion-K\"ahler} if its holonomy subgroup is
contained in
${\rm Sp}(r){\rm Sp}(1)\i SO(4r)$. As usual our manifolds are assumed to be compact and
simply-connected; since ${\rm Sp}(1){\rm Sp}(1)=SO(4)$ we always suppose $r\ge 2$. 
\ind
Despite the terminology, which is unfortunate but classical, a quaternion-K\"ahler
manifold has no natural complex structure: the group ${\rm Sp}(r){\rm
Sp}(1)$ is not contained in $U(2n)$. 
\ind 
 The complex structures $I,J,K$ are not invariant under ${\rm Sp}(1)$, and therefore they
do not correspond any more to parallel complex structures. What remains invariant,
however, is the 3-dimensional space  spanned by $I,J$ and $K$; it gives rise to a 
3-dimensional parallel sub-bundle $E\i{\cal E}nd(T(M))$. The unit sphere bundle $Z\i
E$ is called the {\it twistor space} of $M$; the fibre of $p:Z\rightarrow M$ at a
point $m\in M$ is a sphere ${\bf S}^2$ of complex structures on $T_m(M)$, as in
the hyperk\"ahler case. The link between quaternion-K\"ahler manifolds and
algebraic geometry is provided by the following result of Salamon [S]:
\smallskip 
{\bf Proposition}$.-$ {\it $Z$ admits a natural complex structure, for which the
fibres of $p$ are complex rational curves.}
\ind The construction of this complex structure is quite natural.  Since $E$ is 
parallel, it inherits from the Levi-Civita connection on $T(M)$ a linear
connection, which  is compatible with the metric. It follows that the corresponding
horizontal distribution (App.\ A) induces a horizontal distribution on the fibration
$p:Z\rightarrow M$, that is a sub-bundle
$H\i T(Z)$ which is supplementary to the vertical tangent bundle $T(Z/M)$.
\ind Let $z\in Z$, and let $m=p(z)$. The fibre $p^{-1} (m)$ is canonically
isomorphic to the standard sphere ${\bf S}^2$, and therefore the vertical tangent space
$T_z(Z/M)$ has a well-defined complex structure. The space $H_z$ projects
isomorphically onto
$T_m(M)$, on which $z$ defines by definition a complex structure. The direct
sum of these complex structures define a complex structure on
$T(Z)=T(Z/M)\oplus H$. A non-trivial calculation shows that it is
integrable.\cqfd \smallskip 
\ind As an example, for the quaternionic projective space $M={\bf HP}^r$, the
twistor space $Z$ is ${\bf CP}^{2r+1}$; the fibration $p:Z\rightarrow M$ is 
the natural quotient map $W/{\bf C}^*\rightarrow W/{\bf H}^*$, with $W={\bf
C}^{2r+2}\moins\{0\}={\bf H}^{r+1}\moins\{0\}$.  Its fibres are (complex
projective) lines in ${\bf CP}^{2r+1}$. \smallskip 
\ind The behaviour of the complex manifold $Z$ depends heavily on the
sign of the scalar curvature $k$ of
$(M,g)$. This is a constant; in fact, Berger proved that a $n$\tx dimensional
quaternion-K\"ahler manifold  $(M,g)$ satisfies the Einstein condition
  ${\rm Ric}_g={k\over n}g$ (I refer to [Be], Ch.\ 14.D for a discussion of
the proof). The
case $k=0$ gives the hyperk\"ahler manifolds (\S 3). In the case $k<0$
there seems to be no natural K\"ahler structure on $Z$; actually no compact
example is known. We will therefore concentrate on the case $k>0$,
where some nice geometry appears. Let me recall that a (compact) manifold
$X$ is {\it Fano} if its anticanonical bundle $K_X^{-1}$ is ample (App.\ B).
We will call a quaternion-K\"ahler manifold {\it positive} if its scalar
curvature is positive.
\smallskip 
{\bf Proposition}$.-$ {\it If  $M$ is positive, $Z$ is a Fano manifold
and admits a K\"ahler-Einstein metric}.
\ind The metric on $Z$ is obtained in the same way as the complex structure, by
 putting together the standard metric of the sphere ${\bf S}^2$ on $T(Z/M)$
and the metric of $M$ on $H$ (with the appropriate normalization).\cqfd
\smallskip 
\ind The space $Z$ has one more property, namely a (holomorphic) {\it contact
structure}. We will now explain what this is.
\bigskip
{\bf 4.2. Contact structures}
\smallskip
\ind Let $X$ be a complex manifold. A {\it contact structure} on $X$ is a
 corank 1 sub-bundle $H$ of the (holomorphic) tangent bundle $T(X)$, so that the we
have an exact sequence 
$$0\rightarrow H\longrightarrow T(X)\qfl{\theta } L\rightarrow 0\ ,$$
where $L$ is a line bundle. Moreover 
the 
following equivalent properties must hold:
\indp {\it a}) The 2-form $d\theta $, restricted to $H$, is non-degenerate at
each point\note{The form $d\theta $ is defined locally using a
trivialization of $L$;  it is an easy exercise to check that conditions {\it
a}) and {\it b}) do not depend on the choice of the trivialization.};
\indp {\it b}) $\dim(X)$ is odd, say $=2r+1$, and the form $\theta \wedge
(d\theta) ^r$ is everywhere
$\not=0$;
\indp {\it c}) The $L$\tx valued alternate form $(U,V)\mapsto\theta
([U,V])$ on $H$ is non-degenerate at each point.
\def\et{^{\scriptscriptstyle\times }}
\ind Let  
$L\et$ be the  complement of the zero section in $L^*$. The pull-back of the
line bundle $L$ to $L\et$ has a canonical trivialization, so
$p^*\theta$ becomes a honest 1-form  on
$L\et$. Put $\omega=d(p^*\theta)$. This 2-form is {\it 
equivariant} with respect to the natural action of ${\bf C}^*$
on $L\et$ by homotheties, that is 
$\lambda ^*\omega =\lambda \omega $ for every $\lambda \in{\bf C}^*$.
\smallskip
{\bf Proposition} (``contactization")$.-$ {\it The $2$\tx form $\omega $
 is a  symplectic structure on $L\et$. Conversely, any ${\bf
C}^*$\tx equivariant symplectic $2$\tx form on $L\et$ defines a unique
contact form
$\theta \in H^0(X, \Omega ^1_X\otimes L)$ such that} $\omega =d(p^*\theta)$.
\ind The form $\omega $ is closed, and using {\it b}) above we see easily
that  it is non-degenerate. For the converse, consider the ``Euler field"
$\xi$ on $L\et$ corresponding to the ${\bf C}^*$\tx action. The 1-form 
$i(\xi )\omega $ vanishes on $\xi$ and is equivariant, therefore it is the
pull-back of a form 
$\theta \in H^0(X, \Omega ^1_X\otimes L)$. Since $\omega $ is equivariant, its
Lie derivative $L_\xi \omega $ equals $\omega $; using the Cartan formula
$L_\xi =di(\xi )+i(\xi )d$ we find $\omega =d(p^*\theta )$. It is then
an easy exercise to prove that $\theta $ is a contact form, using for instance
condition {\it a}).\cqfd
 
\smallskip

{\it Example}$.-$ Let $M$ be a complex manifold, and $X={\bf P}T^*(M)$
its (holomorphic) projective cotangent bundle. Recall that the cotangent
bundle $T^*(M)$ has a canonical symplectic structure $\omega=d\eta$,
where $\eta$ is the tautological 1-form on $T^*(M)$: the value of $\eta$
at a point $(\alpha,m)$ of $T^*(M)$ $(m\in M,\alpha\in T_m^*(M))$ is the
pull-back of $\alpha$ by the projection $T^*(M)\rightarrow M$. By
construction $\eta$ is equivariant with respect to the action of ${\bf C}^*$
on $T^*(M)$ by homotheties, and so is $\omega$. By the proposition we
see that $\eta$ is the pull-back of a contact form on $X$. 
\bigskip
\ind Going back to quaternion-K\"ahler manifolds, the link with contact structures
is provided by the following theorem. Part a) is due to Salamon [S], part b) to
LeBrun [L].
\smallskip
{\bf Theorem} (LeBrun, Salamon)$.-$ a) {\it The twistor space of a positive 
quaternion-K\"ahler manifold is a Fano contact manifold, admitting a
K\"ahler-Einstein metric.}
\ind b) {\it Conversely, a Fano contact manifold which admits a
K\"ahler-Einstein metric is the twistor space of a positive 
quaternion-K\"ahler manifold}.
\ind The key point is that the horizontal sub-bundle H$\i T(Z)$ (4.1) is {\it
holomorphic}; this is proved by a local computation, and so is the fact that $H$
defines a contact structure.
\smallskip \ind Thus the classification of positive 
quaternion-K\"ahler manifolds is essentially reduced to a problem of Algebraic
Geometry. We are now going to explain a conjecture describing  this
classification.
\bigskip
{\bf 4.3. Homogeneous contact manifolds}
\smallskip

\ind We  have already mentioned that the only known examples of positive
quaternion-K\"ahler manifolds are symmetric. More precisely, for each simple
compact Lie group $K$ there exists a unique quaternion-K\"ahler symmetric
quotient of $K$; the corresponding
twistor space is homogeneous under the complexification $G$ of $K$. These spaces have
been classified by Wolf [W]. The twistor spaces admit the following
simple description: 
\smallskip  {\bf Proposition}$.-$ {\it Let $G$ be a complex simple Lie group,
${\goth g}$ its Lie algebra. There is a unique closed orbit $X_{\goth g}$ for
the adjoint action of $G$ on ${\bf P}({\goth g})$; $X_{\goth g}$ is a Fano
manifold, and admits a $G$\tx invariant contact structure.}
\smallskip \ind Note that the closure in ${\bf P}({\goth g})$ of any adjoint
orbit contains a closed orbit, necessarily equal to
$X_{\goth g}$. Hence $X_{\goth g}$ is the smallest orbit in ${\bf P}({\goth
g})$.
\smallskip 
{\it Proof}: I will give the proof because it is quite simple, though it
requires some knowledge of algebraic groups.
 Let $X$ be a closed orbit in ${\bf P}({\goth g})$, and let $v$ be a
vector of ${\goth g}$ whose class $[v]\in {\bf P}({\goth g})$ belongs to
$X$.  Since $X$ is projective, the  stabilizer $P$ of $[v]$ contains a Borel
subgroup $B$ of $G$; this means that $v$ is an eigenvector of $B$ in
${\goth g}$. Since ${\goth g}$ is simple, the adjoint representation of $G$
in ${\goth g}$ is irreducible, so $B$ has exactly, up to a scalar, one
eigenvector (``highest weight vector") $v_B\in{\goth g}$; thus $X$ is the $G$\tx
orbit of $[v_B]$. It does not depend on the  particular choice of $B$ because
 all  Borel subgroups are conjugate. 
\ind The pull-back of $X_{\goth g}$ in ${\goth
g}\moins\{0\}$ is an adjoint orbit of $G$; using the Killing form we can view
it as a coadjoint orbit in ${\goth g}^*$. Every such orbit admits a symplectic
form, the Kostant-Kirillov structure, which is ${\bf C}^*$\tx equivariant and $G$\tx
invariant. Using  contactization  we see that $X_{\goth g}$ carries a $G$\tx
invariant contact structure.\cqfd\medskip
\ind  For classical Lie algebras, the contact manifold $X_{\goth g}$ and the
corresponding quaternion-K\"ahler manifold $M_{\goth g}$ are given below:
\smallskip 
$$\hss\vbox{\offinterlineskip\nospacedmath \halign{\tv$#$&\hq$#$\hq\tv&  \hq $#$\hq
\tv  &\hq $#$\hq &\tv#\cr &{\goth g}&X_{\goth g} & M_{\goth g}&\cr\n
&{\goth s}{\goth l}(n) & {\bf P}T^*({\bf P}^{n-1}) & {\bf G}(2, {\bf C}^{n}) &\cr\n
&{\goth o}(n)&{\bf G}_{iso}(2,{\bf C}^n)& 
{\bf G}^+(4,{\bf R}^n)&\cr\n
&{\goth s}{\goth p}(2n) & {\bf CP}^{2n-1} & {\bf G}(1, {\bf H}^{n})={\bf HP}^{n-1} &\cr\n
}}\hss$$ 
\ind We have  described the map
$X_{{\goth s}{\goth p}(2n)} \rightarrow  M_{{\goth s}{\goth p}(2n)}$ in 4.1.
$X_{{\goth o}(n)}$ is the grassmannian of isotropic 2-planes in ${\bf C}^n$ and $M_{{\goth
o}(n)}$ the grassmannian of oriented 4-planes in ${\bf R}^n$; the map $X_{{\goth o}(n)}
\rightarrow  M_{{\goth o}(n)}$ associates to a 2-plane $P\i {\bf C}^n$ the real part of
$P\oplus
\overline{P}$. As in
3.4 we view $X_{{\goth s}{\goth l}(n)}= {\bf P}T^*({\bf P}^{n-1})$ as the space of
flags $D\i H\i {\bf C}^n$, where $D$ is a line and $H$ a hyperplane; choosing a
hermitian scalar product on ${\bf C}^n$, this is also the space
of pairs of orthogonal lines in ${\bf C}^n$. The map $X_{{\goth s}{\goth l}(n)}
\rightarrow  M_{{\goth s}{\goth l}(n)}$ associates to such a pair the 2-plane that they
span.\bigskip
\ind In view of the LeBrun-Salamon theorem (4.2), every
positive quaternion-K\"ahler compact will be symmetric if every Fano contact 
manifold admitting a K\"ahler-Einstein metric is homogeneous. It is tempting to be a
little bit more optimistic and to conjecture:\smallskip 
\indp(C) {\it Every Fano contact manifold is homogeneous}.
\smallskip 
\ind We will give some (weak) evidence for the conjecture. Let $X$ be a  compact
complex manifold, of
dimension
$2r+1$, with a contact structure
$$0\rightarrow H\longrightarrow T(X)\qfl{\theta} L\rightarrow 0\ .$$
 The form
$\theta\wedge(d\theta)^r$ defines a nowhere vanishing section of $K_X\otimes
L^{r+1}$; therefore we have
$K_X\cong L^{-r-1}$, and $X$ is Fano if and only if $L$ is ample.
\smallskip 
{\bf Proposition}$.-$ {\it Let $X$ be a Fano contact manifold. If the 
line bundle $L$ is} very {\it ample, $X$ is homogeneous, and more
precisely  isomorphc to $X_{\goth g}$ for some simple Lie algebra
${\goth g}$}.
\ind {\it Proof}: Let $G$ be the group of automorphisms of $X$ preserving the
contact structure; its Lie algebra ${\goth g}$ consists of the vector fields $V$
on $X$ such that $[V,H]\i H$. Let us  prove that the space of global vector fields
$H^0(X,T(X))$ is the direct sum of ${\goth g}$ and $H^0(X,H)$.
Let $V$ be a vector field on $X$. The map $W\mapsto \theta([V,W])$ from $H$ to
$L$ is ${\cal O}_X$\tx linear, hence by
property {\it c}) of contact structures (4.2),
there exists a unique vector field $V'$ in $H$ such
that $\theta([V,W])=\theta([V',W])$ for all $W$ in $H$. This means that
$[V-V',W]$ belongs to $H$, that is that $V-V'$ belongs to ${\goth g}$. Writing
$V=V'+(V-V')$ provides the required direct sum decomposition.
\ind The map $V\mapsto V'$ provides a ${\bf C}$\tx linear retraction
  of the inclusion of sheaves $H\mono T(M)$; therefore the exact sequence
$$0\rightarrow H\longrightarrow T(X)\qfl{\theta } L\rightarrow 0$$
splits as a sequence of sheaves of vector spaces ({\it not} of ${\cal O}_X$\tx
modules). In particular, the sequence
$$0\rightarrow H^0(X,H)\longrightarrow H^0(X,T(X))\qfl{\theta }
H^0(X,L)\rightarrow 0$$is exact, and $\theta$ {\it induces an isomorphism of
${\goth g}$ onto} $H^0(X,L)$. This isomorphism is equivariant with respect to the
action of $G$. 
\ind We will therefore identify $H^0(X,L)$ with ${\goth g}$.  The diagram of
App.\ B becomes:
$$\xymatrix  @=2pc @C=2.4pc @R=2.4pc  
{ L\et\ar[d]^p \ar[r]^(.4){\mu} &{\goth g}^*\ar@{-->}[d]\\
X\ar@{-->}[r]^(.4){\varphi}&\ {\bf P}({\goth g}^*)\ .
}$$
Let $V\in {\goth g}$. The action of $G$  on $L$ defines a canonical lift
$\widetilde{V}$ of the vector field $V$ to $L\et$. By construction we have
$\langle\mu,V\rangle=\eta(\widetilde{V})$, where $\eta$ is the 1-form
$p^*\theta$ on $L\et$ (4.2). Since 
$\eta $ is preserved by $G$, the Lie derivative
$L_{\widetilde{V}}\eta $ vanishes. By the Cartan homotopy
formula, this implies $$\langle d\mu ,V\rangle =d(i(\widetilde{V})\eta)= 
-i(\widetilde{V})\omega\ ,$$where $\omega:=d\eta$ is the symplectic form on
$L\et$ (this relation means by definition that $\mu$ is a {\it moment map} 
for the action of $G$ on the symplectic manifold $L\et$).  
\ind For $\xi\in L\et$,  $v\in T_\xi(L\et)$, this
 formula reads $\langle T_\xi(\mu)\cdot
v,V\rangle=\omega(v,\widetilde{V}(\xi))$. When $V$ runs in ${\goth g}$, the
vectors $\widetilde{V}(\xi)$ span the tangent space to the orbit $G\xi$ at
$\xi$; thus the kernel of $T_\xi(\mu)$ is the orthogonal of $T_\xi(G\xi)$
with respect to $\omega$. In particular,  if $T_\xi(\mu)$ is injective, the
orbit $G\xi$ is open, and therefore {\it the orbit of $x=p(\xi)$ is open in}
$X$.
\ind Now if $L$ is very ample, $\mu$ is an embedding, hence all the orbits
of $G$ are open -- this is  possible only if $G$ acts transitively on $X$.
Since $X$ is projective this implies that $G$ is semi-simple, so we can
identify ${\goth g}^*$ with ${\goth g}$, and $\varphi(X)$ with a closed
adjoint orbit in ${\bf P}({\goth g})$. It follows easily that ${\goth g}$ is
simple and $\varphi(X)=X_{\goth g}$.\cqfd \medskip
\ind This result is improved in [B3], at the cost of assuming the Lie
algebra ${\goth g}$ {\it reductive} -- this is not too serious since it is
always the case if $X$ admits a K\"ahler-Einstein metric. The main result of
[B3] is:\smallskip  {\bf Theorem}$.-$ {\it Let $X$ be a Fano contact
manifold, such that:
\indp{\rm a)} The rational map $\varphi^{}_L:X\dasharrow {\bf 
P}(H^0(X,L)^*)$  is generically finite {\rm
(}that is, $\dim
\varphi^{} _L(X)=\dim X${\rm );}
\indp{\rm b)} The Lie algebra ${\goth g}$ of infinitesimal contact
automorphisms of
$X$ is reductive.
\ind Then  ${\goth g}$  is simple, and $X$ is isomorphic to $X_{\goth g}$.}
\ind {\it Idea of the proof}: In view of the   above proof, a) implies that
$G$ has an open
 orbit in $L\et$. The image of
this orbit in ${\goth g}$ (identified with ${\goth g}^*$ thanks to b)) is
invariant by homotheties; this implies that it is a nilpotent orbit 
(if a matrix $N$  is conjugate to $\lambda N$ for every $\lambda\in {\bf C}^*$, we
have $\Tr N^p=0$ for each $p$, so $N$ is nilpotent). Thus {\it the image of
$\varphi$ is the closure of a nilpotent orbit in} ${\bf P}({\goth g})$. Then
a detailed study of nilpotent orbits leads to the result.\cqfd

\bigskip
{\bf 4.4. Further developments}
\smallskip
\ind More generally, we can ask which projective varieties admit contact
structures. We have seen two examples, the projective cotangent bundles ${\bf
P}T^*(M)$ (4.2) and the homogeneous spaces $X_{\goth g}$ (4.3). A striking 
fact is that {\it no other example is known}. This leads naturally to the
following question:\smallskip 
 {\it Is every contact manifold isomorphic to a projective cotangent bundle 
or to  one of the homogeneous spaces $X_{\goth g}$}?\smallskip \ind 
This may look overoptimistic, but  let me mention that the answer is positive for:
\indp -- Contact manifolds of dimension $\le 5$: this is due to Ye in
dimension 3 [Ye] and Druel in dimension 5 [D1].
\indp -- Contact toric manifolds [D2]. Druel proves that every such manifold is
isomorphic to  ${\bf P}T^*({\bf
P}^1\times\ldots\times {\bf P}^1)$.
\ind The proofs rely heavily on Mori theory; for this reason they seem
difficult to extend at this point, since Mori theory is well understood only
in low dimension or for toric varieties.

\vfill\eject
\centerline{\bf Appendix A}
\centerline{Connections}\bigskip

\ind Let $M$ be a differentiable manifold, $E$ a vector bundle on $M$, ${\cal
D}iff^1(E)$ the vector bundles of differential operators of order $\le 1$ on 
$E$. A {\it connection} on $E$ is a linear map $\nabla^{}:T(M)\rightarrow
{\cal D}iff^1(E)$ which satisfies the Leibnitz rule
$$\nabla^{}_V(fs)=f\nabla^{}_V(s)+(Vf)s$$
for any vector field $V$, function $f$ and section $s$ of $E$ defined over 
some open subset of $M$. 
\ind The connection extends naturally to the various tensor, symmetric or 
exterior powers of $E$, covariant or contravariant.  For instance, if $b$ is
a bilinear form on
$E$ and $u$ an endomorphism of $E$, we have
$$\nospacedmath\displaylines{
\nabla^{}_V(b)(s,t)=Vb(s,t)-b(\nabla^{}_Vs,t)-b(s,\nabla^{}_Vt)\cr
\nabla^{}_V(u)(s)=\nabla^{}_V(u(s))-u(\nabla^{}_Vs)\cr
}$$for any local sections $s,t$ of $E$. We say that a section $s$ of $E$ (or 
of one of its associated tensor bundles) is {\it parallel} if
$\nabla^{}_Vs=0$ for any vector field
$V$ on $M$.
\smallskip \ind Let $f:M'\rightarrow M$ be a differentiable map. There
exists a natural connection
$f^*\nabla^{}$ on $f^*E$, characterized by the  condition
$(f^*\nabla)^{}_{V'}(f^*s)=f^*(\nabla^{}_Vs)$ for any  section $s$ of $E$ and 
vector fields $V$ on $M$, $V'$ on $M'$ such that $f$ projects $V'$ onto $V$.
In particular, for any path $\gamma:[0,1]\rightarrow M$, we  get a connection
on $\gamma^*E$, or equivalently a first order differential operator
$\nabla^{}_{d/dt}$ of
$\gamma^*E$. Let $p=\gamma(0)$ and $q=\gamma(1)$; given a vector $v_p\in
E_p$, there exists a unique section 
$t\mapsto v(t)$ of $\gamma^*E$  such that
$\nabla^{}_{d/dt}v(t)=0$ and
$v(0)=v_p$. The map
$v_p\mapsto v(1)$ defines the {\it parallel transport} isomorphism
$\varphi_\gamma:E_p\rightarrow E_q$. Observe that a section
$s$ of
$E$ is parallel if and only if $\varphi_\gamma(s(p))=s(q)$ for every  path
$\gamma$ (this implies $s(\gamma(t))=v(t)$, hence $\nabla^{}_{\dot \gamma(t)}s=0$).
\ind The tangent vector $\dot v(0)\in T_{v_p}(E_p)$ is said to be {\it horizontal};
it is easy to show that the horizontal vectors form a sub-bundle $H$ of $T(M)$, the
{\it horizontal distribution} of $\nabla$, which is a supplement of the
vertical sub-bundle $T(E/M)$. 
\medskip
\ind Suppose now $E=T(M)$. The connection is said to be {\it symmetric} (or
torsion-free) if
$\nabla^{}_VW-\nabla^{}_WV=[V,W]$ for any vector fields $V,W$ on $M$. Let $g$ be a
Riemannian metric on $M$; a simple-minded computation shows that {\it there exists a
unique symmetric connection $\nabla^{}$ on $T(M)$ for which $g$ is parallel}. It
is called the Levi-Civita connection of $(M,g)$. 

\vfill\eject
\centerline{\bf Appendix B}
\centerline{Ample line bundles, Hodge theory}\bigskip

{\it Ample line bundles}
\ind  Let $X$ be
a compact complex manifold and $L$ a line bundle on $X$; we suppose
$H^0(X,L)\not= 0$. For
$x\in X$, let
$\varphi_L(x)$ denote the subspace of global sections of $L$ which vanish
at $x$. It is either equal to $H^0(X,L)$ or to a hyperplane in $H^0(X,L)$. In
the first case $x$ belongs to the {\it base locus} $B_L$ of $L$, that is the
subvariety of the common zeros of all sections of $L$. The map $x\mapsto
\varphi_L(x)$ defines a morphism $X\moins B_L\rightarrow {\bf
P}(H^0(X,L))^*$, which we consider as a rational map $X\dasharrow {\bf
P}(H^0(X,L))^*$. We say that $L$ is {\it very ample} if $\varphi_L$ is an
embedding  (this implies in particular $B_L=\varnothing$); it amounts to say
that there is an  embedding of $X$ into some projecive space ${\bf P}$
such that $L$ is the restriction of the tautological line bundle ${\cal
O}_{\bf P}(1)$. We say that
$L$ is {\it ample} if some (positive) power of $L$ is very ample. 
\ind Consider the dual line bundle $p:L^*\rightarrow X$. To any $\xi\in
L^*$  associate the linear form $\mu(\xi): s\mapsto \langle
s(p(\xi)),\xi\rangle$ on $H^0(X,L)$. We have a commutative diagram
$$\xymatrix @=2pc @C=2.4pc @R=2.4pc 
{ L^*\ar[d]^p \ar[r]^(.4){\mu^{}_L} &\ H^0(X,L)^*\ar@{-->}[d]\\
X\ar@{-->}[r]^(.4){\varphi^{}_L}&\ {\bf P}(H^0(X,L)^*)\ .
}$$
 
{\it Hodge decomposition}
\ind Let $X$ be a compact K\"ahler manifold. Recall that a differentiable
form on $X$ is of type $(p,q)$ if it can be written in any system of local
coordinates $(z_1,\ldots,z_n)$ as a sum of forms $a(z,\bar z)dz_{i_1}\wedge
\ldots\wedge dz_{i_p}\wedge  d\bar z_{j_1} \wedge\ldots\wedge d\bar z_{j_q}$.
We denote by $H^{p,q}\i H^{p+q}(X,{\bf C})$ the subspace
 of De Rham cohomology classes of forms of type $(p,q)$; we have
$H^{q,p}=\overline{H^{p,q}}$. The fundamental result of Hodge theory is the
Hodge decomposition
$$H^{n}(X,{\bf C})=\bigoplus_{p+q=n}H^{p,q}\ ,$$
together with the canonical isomorphisms $H^{p,q}\iso H^q(X,\Omega^p_X)$.
In particular,
$$H^{2}(X,{\bf C})=H^{2,0}\oplus H^{1,1}\oplus H^{0,2}\ ,$$with
$H^{2,0}\cong H^0(X,\Omega^2_X)$, embedded into $H^{2}(X,{\bf C})$ by
associating to a holomorphic form its De Rham class.
\ind To any hermitian metric $g$ on $X$ is associated a real
2-form $\omega$ of type $(1,1)$, the {\it K\"ahler form}, defined by
$\omega(V,W)=g(V,JW)$ for any real vector fields $V,W$; the metric is
K\"ahler if $\omega$ is closed. Then its class in $H^{2}(X,{\bf
C})$ is called a K\"ahler class. The K\"ahler classes form an open cone in 
$H^{1,1}_{\bf R}:=H^{1,1}\cap H^{2}(X,{\bf R})$.\ind Let
$L$ be a line bundle on
$X$. The Chern class
$c_1(L)\in H^2(X,{\bf C})$ is integral, that is comes from $H^2(X,{\bf Z})$,
and belongs to $H^{1,1}$. Conversely, any integral class in $H^{1,1}$ is
the Chern class of some line bundle on $X$ (Lefschetz theorem).
\ind  If
$L$ is very ample, its Chern class is the pull-back by
$\varphi_L$ of the Chern class of
${\cal O}_{\bf P}(1)$, which is a K\"ahler class, and therefore $c_1(L)$ is
a K\"ahler class. More generally, if $L$ is ample, some multiple of $c_1(L)$
is a K\"ahler class, hence also $c_1(L)$. Conversely, the celebrated 
Kodaira embedding theorem   asserts that {\it a line bundle whose Chern
class is K\"ahler is ample}. As a corollary, we see that {\it any compact
K\"ahler manifold $X$ with $H^{0}(X,\Omega^2_X)=0$ is projective}: we have
$H^{2}(X,{\bf C})=H^{1,1}$, hence
the cone of K\"ahler classes is open in $H^2(X,{\bf R})$. Therefore it
contains integral classes; by the above results such a class is the first
Chern class of an ample line bundle, hence $X$ is projective. More generally,
the same argument shows that $X$ is projective whenever the subspace $H^{1,1}$
of $H^2(X,{\bf C})$ is defined over ${\bf Q}$.
  \vskip2cm \centerline{ 
REFERENCES} \vglue15pt\baselineskip12.8pt \def\num#1{
\smallskip \item{\hbox to\parindent{\enskip
[#1]\hfill}}} \parindent=1.3cm
\num{A} M.\ {\pc ATIYAH}: {\sl On analytic surfaces with double
points}. Proc.\  R.\  Soc.\  London  {\bf 245}, 237--244 (1958).
\num{Ba} V.\ {\pc BATYREV}: {\sl On the Betti numbers of birationally 
isomorphic  projective varieties with trivial canonical bundles}. Preprint
alg-geom/9710020, to appear in Proc.\ European Algebraic Geometry Conference
(Warwick, 1996).
\num{B1} A.\ {\pc BEAUVILLE}: {\sl 	Vari\'et\'es k\"ahl\'eriennes dont la
premi\`ere classe de Chern est nulle.}  J.\  of Diff.\  Geo\-metry {\bf 18},
755-782  (1983).  
 \num{B2} A.\ {\pc BEAUVILLE}: {\sl   Surfaces $K3$}.  Exp.\ 609
du S\'eminaire Bourbaki; Ast\'erisque {\bf 105}--{\bf 106}, 217--229 (1983). 
\num{B3} A.\ {\pc BEAUVILLE}: {\sl  Fano contact manifolds and nilpotent
orbits}.  Comment.\  Math.\  Helv.\  {\bf 73}, 566--583 (1998).
\num{B-D} A.\ {\pc BEAUVILLE}, R.\ {\pc DONAGI}: {\sl 	La vari\'et\'e des
droites d'une hypersurface cubique de dimension $4$}. C.R.\ 
Acad.\  Sc.\  Paris {\bf 301} (s\'er.\ I), 703--706 (1985). 
\num{Be} A.\ {\pc BESSE}: {\sl Einstein manifolds}. Ergebnisse der
Math.\  {\bf 10}, Springer-Verlag (1987).
\num{B-G} R.\ {\pc BROWN}, A.\ {\pc GRAY}: {\sl Riemannian manifolds
with holonomy group ${\rm Spin}(9)$}. Differential Geometry, in honor
of K.\ Yano, 41--59; Kinokuniya, Tokyo (1972).
\num{B-Y} S.\ {\pc BOCHNER}, K.\ {\pc YANO}: {\sl Curvature and Betti
numbers}. Annals of Math.\  Studies {\bf 32}, Princeton University Press (1953).

\num{Bo1} F.\ {\pc BOGOMOLOV}: {\sl Hamiltonian K\"ahler manifolds}. Soviet
Math.\  Dokl.\  {\bf 19}, 1462--1465 (1978). 
\num{Bo2} F.\ {\pc BOGOMOLOV}: {\sl On the cohomology ring of a simple
hyper-K\"ahler manifold (on the results of Verbitsky)}. Geom.\ Funct.\ Anal.\
{\bf 6}, 612--618   (1996). 
\num{De} O.\ {\pc DEBARRE}: {\sl Un contre-exemple au th\'eor\`eme de Torelli pour 
les vari\'et\'es symplectiques irr\'eductibles}. C.\ R.\ Acad.\ Sci.\ Paris,
  {\bf 299} (s\'er.\ I), 681--684 (1984).
\num{D1} S.\ {\pc DRUEL}: {\sl Structures de contact sur les vari\'et\'es
alg\'ebriques de dimension $5$}. C.R.\  Acad.\  Sci.\  Paris {\bf 327}
(s\'er.\ I),  365--368 (1998).
\num{D2} S.\ {\pc DRUEL}: {\sl Structures de contact sur les vari\'et\'es toriques}.
Math.\  Annalen, to appear (preprint  math/9807030).
\num{F} R.\ {\pc FRIEDMAN}: {\sl  On threefolds with trivial canonical bundle}.
Complex Geometry and Lie theory, 103--134,
Proc.\  Sympos.\  Pure Math.\  {\bf 53}, AMS (1991).
\num{H1} D.\ {\pc HUYBRECHTS}: {\sl Compact hyperk\"ahler
manifolds: basic results}.  Invent.\ math.\ {\bf 135}, 63--113 
(1999).
\num{H2} D.\  {\pc HUYBRECHTS}: {\sl Quelques r\'esultats sur les
vari\'et\'es symplectiques irr\'e\-ductibles}. Proceedings of the
Conference ``Hirzebruch 70" (Warsaw, 1998); Contemp.\  Math., AMS, to
appear.
\num{H-L} D.\ {\pc HUYBRECHTS}, M.\ {\pc LEHN}: {\sl The geometry of moduli 
spaces of sheaves}. Aspects of Math.\ {\bf 31}, Vieweg (1997). 
\num{J1} D.\ {\pc JOYCE}: {\sl Compact Riemannian
$7$\tx manifolds with holonomy $G_2$}. J.\  Diff.\  Geom.\  {\bf 43}, 291--375 (1996).
\num{J2} D.\ {\pc JOYCE}: {\sl Compact Riemannian $8$\tx manifolds with holonomy
${\rm Spin}(7)$}. Invent.\  math.\  {\bf 123}, 507--552 (1996).
\num{J3} D.\ {\pc JOYCE}: {\sl Compact Manifolds with Exceptional Holonomy}.
Doc.\ Math.\ J.\ DMV Extra Volume ICM II, 361--370  (1998).
\num{K} M.\ {\pc KAPRANOV}: {\sl Rozansky-Witten invariants via
Atiyah classes}. Preprint alg-geom/9704009.
\num{K-N} S.\ {\pc KOBAYASHI},  K.\ {\pc NOMIZU}: {\sl Foundations of
differential geometry}.  Wiley-Interscience Publications (1963).
 \num{L} C.\ {\pc LE}{\pc BRUN}: {\sl Fano manifolds, contact
structures, and quaternionic geometry}. Int.\ J.\ of Math. {\bf 6}, 419-437
(1995).   
\num{M} S.\ {\pc MUKAI}: {\sl Symplectic structure of the moduli
space of  sheaves on an abelian or $K3$ surface}. Invent.\ Math.\
{\bf 77}, 101--116  (1984).
\num{OG1} K.\ O'{\pc GRADY}: {\sl The weight-two Hodge structure of 
moduli  spaces of sheaves on a $K3$ surface}. J.\  Algebraic Geom.\ 
{\bf 6}, 599--644 (1997).
\num{OG2} K.\ O'{\pc GRADY}: {\sl Desingularized moduli spaces of 
sheaves on a} $K3$, I and II. Preprints alg-geom/9708009 and
math/9805099.
\num{P} {\sl 	G\'eom\'etrie des surfaces $K3$:
modules et p\'eriodes} (s\'eminaire Palaiseau 81--82).
Ast\'erisque {\bf 126} (1985). 
\num{R-W}  L.\ {\pc ROZANSKY}, E.\ {\pc WITTEN}: {\sl
Hyper-K\"ahler geometry and invariants of three-manifolds}.
Selecta Math.\  {\bf 3},  401--458  (1997).
\num{S} S.\ {\pc SALAMON}: {\sl Quaternionic K\"ahler manifolds}.
Invent.\ Math.\ {\bf 67}, 143--171 (1982).
\num{S-Y-Z} A.\ {\pc STROMINGER}, S.-T.\ {\pc YAU}, E.\ {\pc ZASLOW}:
{\sl   Mirror symmetry is T-duality}. Nuclear Phys.\ 
{\bf B 479},  243--259  (1996).
\num{T} G.\ {\pc TIAN}: {\sl  
Smoothness of the universal deformation space of compact Calabi-Yau manifolds 
and its Petersson-Weil metric}. Math.\ aspects of string theory, 
629--646;  Adv.\  Ser.\  Math.\  Phys.\  {\bf 1}, 
World Scientific, Singapore (1987). 
\num{Va} J.\ {\pc VAROUCHAS}: {\sl Stabilit\'e de la classe des vari\'et\'es
k\"ahl\'eriennes par certains morphismes propres}.  Invent.\ Math.\ {\bf 77},
117--127 (1984). 
\num{V} C. {\pc VOISIN}: {\sl Sym\'etrie miroir}. Panoramas et Synth\`eses 
{\bf 2}, SMF (1996).
\num{W} J.\ {\pc WOLF}: {\sl Complex homogeneous contact manifolds and quaternionic
symmetric spaces}. J.\ Math.\ Mech.\ {\bf 14},  1033--1047 (1965). 
\num{Y} S.-T.\ {\pc
YAU}: {\sl On the Ricci curvature of a compact K\"ahler manifold and the complex
Monge-Amp\`ere equation} I. Comm.\ Pure and Appl.\ Math.\ {\bf 31},
339--411 (1978).
\num{Ye} Y.-G.\ {\pc YE}: {\sl A note on complex projective threefolds admitting
holomorphic contact structures}. Invent.\ Math.\ {\bf 115}, 311--314  (1994).

 \vskip1cm \def\pc#1{\eightrm#1\sixrm}
\hfill\vtop{\eightrm\hbox to 5cm{\hfill Arnaud {\pc BEAUVILLE}\hfill}
 \hbox to 5cm{\hfill DMI -- \'Ecole Normale
Sup\'erieure\hfill} \hbox to 5cm{\hfill (UMR 8553 du
CNRS)\hfill}
\hbox to 5cm{\hfill  45 rue d'Ulm\hfill}
\hbox to 5cm{\hfill F-75230 {\pc PARIS} Cedex 05\hfill}}
\end